\documentclass[12pt]{article}
\usepackage{amsmath}
\usepackage{amsthm}
\advance \topmargin by -\headheight
\advance \topmargin by -\headsep
\evensidemargin \oddsidemargin
\theoremstyle{plain}
\newtheorem{thm}{Theorem}[section]
\newtheorem{prop}[thm]{Proposition}
\newtheorem{lem}[thm]{Lemma}
\newtheorem{cor}[thm]{Corollary}
\newtheorem{conj}[thm]{Conjecture}
\newtheorem*{mainthm}{Main Theorem}
\theoremstyle{remark}
\newtheorem{rmk}[thm]{Remark}

\theoremstyle{definition}

\title{\Large \bf Homologically Trivial Group Actions on 4-Manifolds}
\author{ Allan L. Edmonds\\
Department of Mathematics\\[-10pt]
Indiana University\\[-10pt]
Bloomington, IN 47405\\[-10pt]
\texttt{edmonds@indiana.edu}}
\date{}

\begin{document}
\maketitle

\section{\large Introduction}In this paper we investigate pseudofree homologically
trivial actions of {\em arbitrary} finite groups on simply connected $4$-manifolds. 
It turns out that there are surprisingly strong restrictions on the groups and on the
actions, especially when one considers the fact that any finite group can act, even
freely, on some simply connected $4$-manifold.

Recall that a finite group action is \emph{free} if no nontrivial group element has a fixed point, \emph{pseudofree} if the
action is free on the complement of a discrete set, and
\emph{semifree} if it is free on the complement of the fixed point set of the whole group.  An action is \emph{locally linear} if
each point has a neighborhood invariant under the isotropy group at the point on which the action is equivalent to a linear
action of the isotropy group on some euclidean space.  An an action is said to be \emph{homologically trivial} if the induced
action on all integral homology groups is the identity.

\begin{mainthm} If a finite group $G$ acts locally linearly, pseudofreely, and
homologically trivially, on a closed, simply connected $4$-manifold $X$, with
$b_2(X)\ge 3$, then $G$ is cyclic and acts semifreely, and the fixed point set consists of $b_2(X)+2$ isolated points.
\end{mainthm} 

When the second Betti number $b_2(X)\le 2$ other groups do arise.  In the case of $S^4$ (corresponding to the case when
$b_2=0$),a number of other groups do come up, since any free linear action on on
$S^3$ can be suspended to a pseudofree action, so a number of other groups arise. 
It turns out every group that acts pseudofreely on $S^4$ arises in this way.  Most
pseudofree actions on the complex projective plane
$\mathbf{C} P^2$ ($b_2=1$) are in fact semifree, as follows from the classification of D.
Wilczy\'nski \cite{Darek87}, \cite{Darek90}  and of I. Hambleton and R. Lee
\cite{Hambleton88}.  But there is a pseudofree, non-semifree action of
$C_3\times C_3$.  Also there are standard pseudofree actions of the polyhedral
groups
$C_n$,
$D_{2n}$, $A_4=T$, $S_4=O$, and $A_5=I$ on
$S^2$.  These non-cyclic groups give rise to pseudofree, non-semifree actions of
these groups on
$S^2\times S^2$  ($b_2=2$).  M. McCooey \cite{McCooey98} has recently analyzed the groups
that act pseudofreely on
$S^2\times S^2$.  It turns out that every  finite group that acts pseudofreely on  $S^2\times S^2$ also admits a linear action on
$\mathbf{R}^6$ leaving the standard
$S^2\times S^2$ invariant, and, of course, these groups can be explicitly catalogued.

Note also that most finite cyclic groups of odd order admit such actions on any
closed simply connected $4$-manifold by results of A. Edmonds \cite{Edmonds87}.  (There are additional
restrictions on even order periodic maps originating in Rokhlin's Theorem.  And
there are certain cases when the order is divisible by 3 in which the group action
must have a two-dimensional fixed point set.)


It remains for future work to consider actions with $2$-dimensional fixed singular
sets, in general.

Here is a summary of the main argument.  We observe that the theorem is true for cyclic
groups.  We therefore attempt to proceed by induction on the order of the group.  This
requires separating the two conclusions (cyclic and semifree) and the consideration of some
minimal non-cyclic groups.  In Theorem
\ref{cyclic} we show that a group acting semifreely and homologically trivially must be cyclic,
provided
$b_2\ge 1$.  Then in Theorem
\ref{elemabel} it is shown that a homologically trivial $C_p\times C_p$ action has nonempty
fixed point set, provided $b_2\ge 2$.  In Theorem \ref{metacyclic} we rule out pseudofree, homologically trivial
actions of nonabelian metacyclic groups when $b_2\ge 3$.  Finally the proof of the Main
Theorem is completed in Theorem \ref{semifree} where it is shown that a group acting
pseudofreely and homologically trivially must in fact semifreely, provided $b_2\ge 3$.  At this
last stage we apply the full induction hypothesis to consider groups in which every maximal
subgroup is cyclic and acts semifreely, to reach the final conclusion, applying Theorems \ref{elemabel} and \ref{metacyclic}.

\section{\large The Borel Fibering and Spectral Sequence} 
We follow Bredon \cite{Bredon72},
Chapter VII, or tom Dieck \cite{tomDieck87}, Chapter III, and refer to those sources for more
details. Let
$G$ be a finite group.  Let
$E_G\to B_G$ be a universal principal $G$-bundle.  We think of $E_G$ as a
contractible space on which
$G$ acts freely (on the right) with orbit space $B_G$.  Although $E_G$ and $B_G$ are
infinite dimensional they can be chosen to have finite $n$-skeleton for all $n$.  The
classifying space $B_G$ is an Eilenberg-MacLane space of type $(G,1)$.

Now let $X$ be a (reasonable) left $G$-space, say of the $G$-homotopy type of a finite $G$-CW complex.  We let $$X_G=E_G
\times_G X$$ This is the bundle with base $B_G$ and fiber $X$ associated to the principal bundle.

One then has the spectral sequence of the fibering:
$$ E_2^{i,j}=H^i(B_G;H^j(X;K))\Rightarrow H^{i+j}(X_G;K)
$$ Here $K$ denotes any abelian coefficient group or commutative ring.  It is
understood that as a coefficient group
$H^j(X_G;K)$ is a module over $K[G]$ and we have cohomology with twisted
coefficients.  The differentials $$d_r:E_r^{i,j}\to E_r^{i+r,j-r+1}$$ have bidegree
$(r,-r+1)$.

The limit term of this spectral sequence $E_\infty^{**}$ is the graded group
associated to some filtration of $H^*(X;K)$:
$$ E_\infty^{i,j}=\frac{\mathcal{F}_jH^{i+j}(X_G;K)}{\mathcal{F}_{j-1}H^{i+j}(X_G;K)}
$$

It is important to emphasize that $\bigoplus_{i+j=k}E_\infty^{i,j}$ is not isomorphic
to $H^k(X_G)$ in general.  But if $K$ is a field, then the two vector spaces have the
same dimension.  And if $K$ is a finite abelian group, then the two groups have the
same cardinality, a fact we shall use.


The $G$-equivariant projection $E_G\times X\to X$ induces a map 
$$\varphi:X_G\to X^*$$ where $X^*$ denotes the orbit space.  
For $x^*=\varphi([e,x])$, one can identify $\varphi^{-1}(x^*)=B_{G_x}$.

Let $F=X^G$ denote the
fixed point set of $G$ on
$X$.  We can regard $F$ as being embedded in both $X$ and $X^*$.  Also
$$X_G\supset F_G = B_G\times F$$  

Let $S$ denote the singular set for the action of $G$ on
$X$ consisting of all points with nontrivial isotropy group.  For any closed invariant
subset $A\subset X$ there is a corresponding relative spectral sequence.   Since the projection 
$\varphi:X_G\to X^*$ has point inverses of the form $\varphi(x^*)=B_{G_x}$, which is acyclic for
$x\not\in S$, one obtains the following result from a suitable form of the Vietoris-Begle Mapping Theorem.

\begin{prop} The induced map $$\varphi^* :H^q(X^*,A^*;K)\to H^q(X_G,A_G;K)$$ is an
isomorphism for all closed invariant $A\subset X$
 containing $S$.
\end{prop}

Since $H^q(X^*)=0$ for $q>\dim X$, we have the following conclusion.

\begin{cor} \label{l:inciso}Inclusion induces an isomorphism $H^q(X_G)\to H^q(S_G)$ for $q$
greater than the dimension of $X$.
\end{cor}

In particular, if $G$ acts semifreely, then  $H^q(X_G)\approx H^q(B_G\times F)$ for $q>\dim(X)$.

\section{\large The Borel Fibering Specialized to Dimension 4}   Now suppose that
$X$ is a closed, simply connected $4$-manifold on which $G$ acts homologically
trivially.  
Then the spectral sequence has at most 3 nonzero rows ($j=0,2,4$) and the only
possibly nontrivial differentials are then
$d_3$ of bidegree $(3,-2)$ and $d_5$ of bidegree $(5,-4)$.  If, in addition, $H^*(G;K)$ vanishes in odd degrees, then the spectral
sequence collapses for trivial formal reasons.    A relevant example, for $K=\mathbf{Z}$, of this is the case of a 4-periodic group
$G$, which necessarily has $H^3(G;\mathbf{Z})=0$.  But for groups like $C_p\times C_p$, the spectral sequence need not always
collapse.

A key point is that if $H^2(X;K)\ne0$, then the spectral sequence collapses if and
only if $d_3:H^2(X;K)\to H^3(G;K)$ vanishes.  (Here we identify $E_2^{02}=E_3^{02}=H^0(G;H^2(X;K))=H^2(X;K)$ and
$E_2^{30}=E_3^{30}=H^3(G;H^0(X;K))=H^3(G;K)$.)  This follows from the multiplicative structure of the spectral sequence, since
the elements of $H^4(X;K)$ are products of elements in $H^2(X;K)$.

Without any further assumptions on $G$ or $X$ one can argue that the spectral sequence collapses for $G$ for $\mathbf{Z}_p$
coefficients, provided $p\ne 2$ or $3$, and $(6,|G|)=1$, but when $p=2$ or $3$,  more care is needed.   See the proof of
Proposition \ref{elemabel}, where we consider the case of $G=C_p\times C_p$ and $K=\mathbf{Z}_p$.

%

%

\section{\large A Little Group Theory} Here we prove a fact about maximal proper
subgroups of a general finite group, useful, as we shall see, for constructing certain
proofs by induction on the order of a group.  Although this result would be well known to finite group theorists we have no
easilly accessible reference and so include a short proof.

\begin{prop} \label{max}Let $G$ be a finite group.  Then either $G$ contains a maximal proper
subgroup that is normal in $G$ or $G$ contains two distinct maximal subgroups that
have a nontrivial intersection.
\end{prop}

\begin{rmk} Of course the second alternative must hold for a simple group.  The
example of a dihedral group $D_{2p}$ shows that the second alternative does not
always hold.
\end{rmk}

\begin{proof} We suppose that no maximal subgroup is normal and that any two
distinct maximal subgroups intersect only in the identity element and derive a
contradiction.

Let $M_{i,j}$, $i=1,\dots, r$, $j=1,\dots, s_i$ be the full list of maximal subgroups of
$G$, indexed so that $M_{i,j}$ and $M_{k,\ell}$ are conjugate if and only if
$i=k$.  We let $m_{i,j}$ denote the order of $M_{i,j}$.  Since distinct maximal
subgroups intersect only in the identity element, we see that $s_i=|G/M_{i,j}|$, that
is $G$ acts transitively by conjugation on
$\{M_{i,j}:j=1,\dots, s_i\}$, with isotropy group $M_{i,j}$.

We now count the total number $|G|$ of elements in $G$:
\begin{align} |G|-1 &= \sum_{i=1}^r\sum_{j=1}^{s_i}m_{i,j}-1 \notag \\ 
    &= \sum_{i=1}^r {s_i}(m_{i,j}-1) \notag \\
    &= \sum_{i=1}^r |G|-s_i \notag
\end{align} where $m_{i,j}=|M_{i,j}|$, which we next abbreviate to $m_i$.  Dividing
through by $|G|$ we obtain
\begin{align} 1-\frac{1}{|G|} &= \sum_{i=1}^r 1-\frac{1}{m_i} \notag 
\end{align} But $2\le m_i < \frac{|G|}{2}$, so each term on the right is at least
$\frac{1}{2}$.  Therefore the right hand side is at least $1$, while the left hand side
is less than
$1$, unless $r=1$.  But, if $r=1$, that is all maximal subgroups are conjugate, we
have
$$ 
1-\frac{1}{|G|}= 1-\frac{1}{m_1}
$$ 
But this would mean that a maximal subgroup has the same order as the group
itself, a contradiction.
\end{proof}

\section{\large Semifree Actions} If a finite group $G$ acts
pseudofreely with a nonempty fixed point set, then it acts freely on a small
invariant 3-sphere around the fixed point.  In particular it follows that such
a group has periodic cohomology.  We will show that such groups cannot act
pseudofreely unless they are cyclic.  We begin with the computation of such a
group's integral cohomology.  See, for example Brown \cite{Brown82}, Section VI.9, including
exercises $3$ and $4$.

\begin{prop}\label{periodic} If $G$ is a finite group of with cohomology of period 4,
then 
$$H^i(G;\mathbf{Z}) =
\begin{cases}
\mathbf{Z} & (i=0)\\  G/[G,G] & (i >0, i=4k+2\\ 
\mathbf{Z}/|G| & (i>0, i=4k) \\ 0 & (\text{otherwise})
\end{cases}$$\qed
\end{prop}

\begin{thm} \label{cyclic} If a finite group $G$ acts semifreely and homologically
trivially  on a simply connected
$4$-manifold $X$ such that $b_2(X)\ge 1$, then $G$ is cyclic.
\end{thm}

\begin{rmk} Any group that acts freely (or semifreely) on $S^3$ acts semifreely on
$S^4$, by suspension.  So, when $b_2=0$, other groups do arise.
\end{rmk}

\begin{proof}  If the fixed point set $X^G$ contains a $2$-dimensional component,
then
$G$ acts effectively, preserving orientation, on the normal $2$-plane to the
$2$-dimensional fixed point set, implying that $G$ is indeed cyclic.  Henceforth we
concentrate on the case when $X^G$ is finite.

For any $C_p<G$, the fixed sets $X^G$ and $X^{C_p}$ are equal, since the action is
semifree, and consist of exactly $b_2(X)+2$ points. 


 Since $G$ has an isolated fixed point, $G$ is a spherical space form group, and in
particular $G$ has periodic cohomology of period $4$ (or $2$ if $G$ is in fact cyclic),
concentrated in even dimensions, as described above.

Now
$H^6(X_G)=H^6(X^G\times B_G)$, and both are finite.  Since the cohomology of $X$
and of $G$ are both concentrated in even dimensions, the Borel spectral sequence
collapses, and we have
\begin{align} |H^6(X_G)| &=|H^6(G)|\times |H^4(G; H^2(X))|\times |H^2(G)|
\notag \\
          &=|G/[G,G]|^2\times |G|^{m} \notag
\end{align} where $m=b_2(X)$.  On the other hand,
\begin{align} |H^6(X^G\times B_G)| &= |H^0(X^G)\otimes H^6(G)| \notag \\  
                    &=|G/[G,G]|^{m+2} \notag
\end{align} where $m+2$ is the number of fixed points, all of which are isolated. 

It follows that
$|G/[G,G]| = |G|$ and hence that $[G,G]$ is trivial.  Thus $G$ is abelian.  But an abelian
periodic group is cyclic.
\end{proof}


\section{\large Elementary Abelian $p$-Groups}
We will apply a result of A. Borel  \cite{Borel60}, Chapter XII, Theorem 3.4:  If an elementary abelian $p$-group acts in a
reasonable manner on a reasonable space $X$ such that the action is trivial on $\bmod\ p$ homology and the $\bmod\ p$
spectral sequence of the fibering $X_G\to B_G$ collapses, then the fixed point set $X^G$ is nonempty and $\dim
H^*(X^G;\mathbf{Z}_p) = \dim H^*(X;\mathbf{Z}_p)$.  See also tom Dieck \cite{tomDieck87}, Chapter III, pp. 200-201.

The following is the result we need about the action of an elementary abelian
$p$-group of rank 2.  

\begin{prop}\label{elemabel} Suppose that the group $G=C_p\times C_p$ of order
$p^2$, where $p$ is prime, acts effectively and homologically trivially on a closed,
simply connected
$4$-manifold
$X$ such that
$m=b_2(X)\ne 1$ if $p=3$, and $m=b_2(X)\ne 0,2$ if $p=2$.  Then the fixed set
$X^G$ consists of $m+2$ isolated points.
\end{prop}

\begin{proof} We fix an elementary abelian
$p$-group $G$ of rank $2$ acting effectively and homologically trivially on a closed,
simply connected
$4$-manifold
$X$.  The fixed point set must be finite.  Otherwise $G$ would act, preserving orientation, on a normal
$2$-disk to the fixed point set, implying that $G$ is cyclic.

We analyze the Borel spectral sequence (with coefficients $\mathbf{Z}_p$) for the
fibering
$$X\to X_G\to B_G$$ We have
$E_2^{i,j}=H^i(G;H^j(X;\mathbf{Z}_p))\Rightarrow H^{i+j}(X_G;\mathbf{Z}_p)$.  According to the result of Borel mentioned above,
it suffices to show that the spectral sequence collapses, that is, 
$E_2=E_\infty$.  For this consideration of the ring sturcture of the spectral sequence implies that it suffices to show that
$d_3:E_3^{0,2}\to E_3^{3,0}$,
$d_3:E_3^{0,4}\to E_3^{3,2}$, and
$d_5:E_5^{0,4}\to E_5^{5,0}$ all vanish.

Let $u\in E_3^{0,2}=H^2(X;\mathbf{Z}_p)$ be nonzero.  Suppose that $u^2\ne 0$.  Of
course
$u^3=0$.  Then
$0=d_3(u^3)=3d_3(u)\cdot u^2$.  It follows that $d_3(u)=0$, provided $p\ne 3$. 
(We'll deal with the case
$p=3$ separately below.)

Suppose $u^2=0$.  Then $0=d_3(u^2)=2d_3(u)\cdot u$, implying that $d_3(u)=0$,
provided
$p\ne 2$.   (We'll also deal with the case $p=2$ separately below.)  This shows that
$d_3$ vanishes on $E_3^{0,2}$.  But
$E_3^{0,4}=H^4(X;\mathbf{Z}_p)$ is generated by elements of degree $2$, provided
$H^2(X)\ne0$.  In any case, if $v\in H^4(X;\mathbf{Z}_p)$,
$0=d_3(v^2)=2d_3(v)\cdot v$, so that $d_3$ vanishes on $E_3^{0,4}$, provided
$p\ne 2$.

Now we know that $E_2=E_5$.  It remains to see that $d_5$ vanishes.  It suffices to
show this for
$d_5:H^4(X;\mathbf{Z}_p)\to H^5(G;\mathbf{Z}_p)$.  But this follows automatically if
$H^2(X)\ne 0$, so that the latter group generates $H^4(X)$.  And if
$H^2(X)=0$, it also follows for $p\ne 2$, as for $d_3$.

Now consider more carefully the case that $p=3$.  We only need to show that $d_3$
vanishes on $E_3^{0,2}$.  So let
$u\in H^2(X;\mathbf{Z}_p)$.  If $u^2=0$, then the previous argument suffices.  So
suppose
$u^2\ne 0\pmod p$.  Notice that, since $b_2(X)>1$, there is necessarily an element
$v\in H^2(X;\mathbf{Z}_p)$ such that $u$ and $v$ are independent and
$uv=0$.  Since $u$ and $v$ are independent over $\mathbf{Z}_p$ in
$H^2(X;\mathbf{Z}_p)$, they are also independent in
$H^2(X_G;\mathbf{Z}_p)$ over $H^*(B_G;\mathbf{Z}_p)$.  But, 
$$0=d_3(uv)=du\cdot v + u\cdot dv$$ shows that $u$ and $v$ are dependent unless
$du=0$, as required.  

Now consider more carefully the case that $p=2$.  To show that $d_3$ vanishes we
only need to show that
$d_3(u)=0$ for nonzero $u\in H^2(X;\mathbf{Z}_2)$ such that $u^2=0 \pmod 2$. 
Clearly then
$b_2(X)\ne 1$, and the hypothesis that $b_2\ne 2$ shows that we can find a
nonzero
$v\in H^2(X;\mathbf{Z}_2)$ independent of $u$ such that
$uv=0$, i.e., $u^\perp\ne\left< u\right>$.  Then, again, $$0=d_3(uv)=du\cdot v +
u\cdot dv$$ which shows that $u$ and $v$ are dependent in
$H_G^2(X;\mathbf{Z}_p)$, over
$H^*(B_G;\mathbf{Z}_p)$, unless
$du=0$, as required. 

Since the $\bmod\ p$ spectral sequence collapses and $G$ acts trivially on $\bmod\ p$
homology, the result of Borel \cite{Borel60} mentioned above implies that the fixed set $X^G\ne \emptyset$ and that
$\dim H^*(X^G;\mathbf{Z}_p)=\dim H^*(X;\mathbf{Z}_p)$.  Since we have already
observed that $X^G$ must be finite, the theorem is proved.
\end{proof}

\begin{rmk} The exceptions when $p=3$ or $p=2$ are necessary.  There is a
$C_3\times C_3$ action on $\mathbf{C} P^2$ with no fixed points.  And there is a
$C_2\times C_2$ action on $S^2$ with no fixed points, leading to a $C_2\times C_2$
action on $S^2\times S^2$ with no fixed points.  Also it is easy to construct a linear
$C_2\times C_2$ action on
$S^4$ with no fixed points.
\end{rmk}


\begin{rmk} An action as above cannot be semifree or pseudofree.  For the existence of an
isolated fixed point then implies that $G$ acts freely on
$S^3$ (at least up to homotopy if the action is not locally linear), so that $G$ would
have periodic cohomology.  But, of course, $C_p\times C_p$ does not have periodic
cohomology.
\end{rmk}

\begin{rmk}
We note that a result that implies Proposition \ref{elemabel} appears in A. Assadi
\cite{Assadi90}, Theorem 5.1.  But the proof there is unconvincing because it
makes no allowance for the counterexamples mentioned above when $b_2\le 2$ and $p=2$ or $3$.
With due care, however, that proof, which proceeds along lines quite different from the proof given here, can be made to work.
\end{rmk}

%

\section{\large Metacyclic Groups}

Next we examine the action of a nonabelian metacyclic group.  In this case we do
assume the action is pseudofree.  If it were true that the entire group had a fixed
point, then most such groups would be easily ruled out, since all but the dihedral groups have
period greater than 4, while a group acting freely on a $3$-sphere around a fixed
point would have to have period 4.  And a dihedral group cannot act freely on $S^3$, by Milnor's $2p$ condition 
\cite{Milnor57}.  We begin with a computation of the cohomology of such a group.

In the rest of this section we suppose that $G$ is a nonabelian extension with
normal subgroup
$C_p$ and quotient group $C_q$, where $p$ and $q$ are (distinct) primes.  Then there is a
short exact sequence

$$1\to C_p \to G \to C_q \to 1$$ and $G$ has presentation of the form
$$\left<  a, b : a^p=1, b^q=a^s, bab^{-1}=a^r\right>$$ where $r^q\equiv 1\bmod p$. 
Since 
$G$ is nonabelian, we have
$r\not\equiv 1 \bmod p$.  Moreover, the extension must in fact be split, i.e., $s=0$. 
For, if
$b^q$ were nontrivial, then $b$ would generate the whole group as a cyclic group,
since we are assuming that $p$ is prime.

We will need the cohomology of $G$ with coefficients in the following $G$-modules: 
$\mathbf{Z}$, $\mathbf{Z}[G]$, $\mathbf{Z}[C_q]=\mathbf{Z}[G/C_p]$, and $\mathbf{Z}[G/C_q]$ (where we choose a
representative non-normal subgroup $C_q< G$), which correspond to the cohomology of an orbit, for each possible orbit $G/H$
of
$G$.  Note in particular that $G$ has periodic cohomology of period $2q$.

\begin{lem} Let $G$ be a nonabelian extension of a normal cyclic subgroup $C_p$ of
prime order $p$ by a cyclic quotient group $C_q$ of prime order $q$.  Then the
cohomology of $G$ is given as follows:
\begin{align}  H^i(G;\mathbf{Z}) &=
\begin{cases}
\mathbf{Z} & (i=0)\\ 
\mathbf{Z}_q & (i >0, i=2k, k\not\equiv 0\bmod q\\ 
\mathbf{Z}_{pq} & (i=2\ell q) \\ 0 & (\text{otherwise})
\end{cases}
\notag\\ H^i(G;\mathbf{Z}[G]) &=
\begin{cases}
\mathbf{Z} & \text{for } i=0 \\ 0 & \text{for } i \ne 0
\end{cases}
\notag\\ H^i(G;\mathbf{Z}[G/C_q])&=
\begin{cases}
\mathbf{Z} & (i=0) \\
\mathbf{Z}_{q} & (i=2k>0) \\ 0 & (\text{otherwise})
\end{cases}
\notag\\ H^i(G;\mathbf{Z}[G/C_p]) &=
\begin{cases}
\mathbf{Z}  & (i=0) \\
\mathbf{Z}_{p} & (i=2k>0) \\ 0 & (\text{otherwise})
\end{cases}
\notag
\end{align}

\end{lem}

\begin{proof}[Proof Sketch] Recall that  the automorphism $C_p\to C_p$ given by
$x\mapsto x^r$ induces multiplication by
$r^k$ on $H^{2k}(C_p)=\mathbf{Z}_p$.   This is best seen by using the ring structure
on
$H^*(C_p;\mathbf{Z}_p)=\Lambda(s_1)\otimes \mathbf{Z}_p[t_2]$, where the
degree $2$ generator $t_2$ is the image of the degree $1$ generator under the
$\bmod p$ Bockstein. Given this observation, these calculations can all be done using
the Lyndon-Hochschild-Serre spectral sequence for the group extension $$1\to C_p
\to G \to C_q \to 1$$ Further details are omitted.
\end{proof}

\begin{thm}\label{metacyclic} Suppose that a nonabelian metacyclic group $G$ of
order $pq$, where $p$ and $q$ are distinct primes, acts pseudofreely and
homologically trivially on a closed, simply connected
$4$-manifold
$X$, then  $b_2(X)\le 2$.
\end{thm}

\begin{proof}

 We analyze the spectral sequence of the Borel fibering $X\to X_G\to B_G$, making
use of  the cohomology groups of $G$ with coefficients in
$H^0(G/H)$ as $H$ varies over the subgroups of $G$.

In particular, we calculate $H^{*}(X_G)$ and compare the results of that calculation
with a direct calculation of $H^{*}(S_G)$, using the fact that these two cohomology
groups agree in dimensions greater than $4$.  where
$G$ acts on
$X$ with singular set
$S$ consisting of a finite union of orbits with isotropy groups $G$,
$C_p\triangleleft G$, or
$C_q<G$.

The spectral sequence for $H^*(X_G)$ has $E_2$-term
$$ E_2^{i,j}=H^i(G; H^j(X))
$$ Since $H^j(X)\ne 0$ only for $j=0,2,4$, and $H^i(G)\ne 0$ only for $i$ even, we see
that the spectral sequence collapses.  Since the action of $G$ on $X$ is homologically
trivial, we find that for any $n>4$ the cohomology group $H^n(X_G)$ has a filtration
with associated graded group $\mathcal{G}H^{n}(X_G)= (H^{n-4}(G), H^{n-2} (G)^{m},
H^n(G))$, where
$m=b_2(X)$. 

First consider the case that $n=4q>4$. Then $\mathcal{G}H^{4q}(X_G)= 
({\mathbf{Z}_{pq}, \mathbf{Z}_{q}}^{m}, \mathbf{Z}_{pq})$. In particular 
$H^{4q}(X_G)$ is finite of order
$p^{2}q^{m+2}$.

On the other hand, we can compute the cohomology of $S_G$ directly.  The singular
set
$S$ is a union of
$x_p$ orbits of size
$p$ (isotropy group conjugate to $C_q$), $x_q$ orbits of size $q$ (isotropy group
$C_p$), and $x_1$ orbits of size 
$1$ (isotropy group $G$).  Each orbit of size $p$ contributes a summand
$H^{4q}(G;\mathbf{Z}[G/C_q])=\mathbf{Z}_q$; each orbit of size
$q$ contributes $H^{4q}(G; \mathbf{Z}[G/C_p]= \mathbf{Z}_p$; and each orbit of size
$1$ contributes a summand
$H^{4q}(G;\mathbf{Z})={\mathbf{Z}_{pq}}$.  It follows that
$H^{4q}(S_G)=\mathbf{Z}_q^{x_p}\oplus
\mathbf{Z}_p^{x_q}\oplus \mathbf{Z}_{pq}^{x_1}$, which has order $(pq)^{x_1}q^{x_p}p^{x_q}=p^{x_1+x_q}q^{x_1+x_p}$.  Since
$H^{4q}(X_G)\approx H^{4q}(S_G)$ by Corollary \ref{l:inciso}, we conclude that  $x_1+x_q=2$ and
$x_1+x_p=m+2$.

Next we consider the case that $n=4q+2>4$. Then $\mathcal{G}H^{4q+2}(X_G)= 
({\mathbf{Z}_{q}, \mathbf{Z}_{pq}}^{m}, \mathbf{Z}_{q})$. In particular 
$H^{4q+2}(X_G)$ is finite of order
$p^{m}q^{m+2}$.  But we may  also compute directly
that $H^{4q+2}(S_G)=\mathbf{Z}_q^{x_1}\oplus
\mathbf{Z}_p^{x_q}\oplus \mathbf{Z}_q^{x_p}$, which has order
$p^{x_q}q^{x_1+x_p}$.  Therefore, $x_q=m$ and $x_1+x_p=m+2$.

These two sets of equations imply immediately that $m=x_q\le x_1+x_q=2$  Since $m=b_2(X)$, we are done.
\end{proof}

\section{\large The Main Theorem}

In this section we finish the proof of the main theorem by showing that
homologically trivial, pseudofree actions are semifree, provided $b_2\ge 3$.  The
argument will proceed by induction on the order of the group.  The result is easy
for cyclic groups, but we also need the analyses of the preceding two sections about
elementary abelian groups and metacyclic groups.

\begin{thm} \label{semifree} Suppose that a finite group $G$ acts pseudofreely and
homologically trivially on a closed, simply connected $4$-manifold $X$ such that $b_2(X)\ge
3$.  Then $G$ acts semifreely.
\end{thm}
\begin{proof} For each nontrivial element $g\in G$, the
fixed point set
$X^g$ consists of exactly $m+2$ points, where $m=b_2(X)$.  To see this, let $\Lambda(g)$ denote the Lefschetz number of $g$.  
Then, by the Lefschetz Fixed Point Formula we have
$$
|X^g|=\chi(X^g)=\Lambda(g)=\chi(X)=b_2(X)+2
$$  
Say that a
subgroup
$H$ of
$G$ is {\em good} provided that its fixed point set $X^H$ consists of exactly $m+2$
points.  In other words $X^H=X^h$ for all nontrivial $h\in H$.  (One always has
$X^H\subset X^h$.)  We need to show that
$G$ itself is good.  We proceed inductively toward this end.  Since every cyclic
subgroup is good, we may assume that every proper subgroup is good.  By
induction we know that such a subgroup $H$ acts semifreely; then by Theorem \ref{cyclic} $H$
must be a cyclic group
$C_k$ acting semifreely on
$X$ with
$m+2$ isolated fixed points, where $m=b_2(X)$. 

\noindent\textbf{ Case 1.} \emph{There are two distinct maximal subgroups that intersect nontrivially. } Suppose that $H_1$
and
$H_2$ are maximal subgroups and that $h\in H_1\cap H_2 \backslash \{ e\}$.  It follows easily that the
subgroup
$H_1H_2$ generated by $H_1$ and $H_2$ is also good, since $H_1$, $H_2$ and $\left<h\right>$ all have the same fixed point
set.   But maximality implies this good subgroup is $G$ itself.

\bigskip
\noindent\textbf{Case 2.} \emph{All maximal proper subgroups intersect trivially. } Then, by Proposition \ref{max} there is a
maximal proper subgroup
$H<G$ that is normal in $G$.  By maximality it follows that
$G/H\approx C_q$, a cyclic group of prime order $q$.  By induction, $H$ acts semifreely, and so by Theorem~\ref{cyclic} $H$ is
cyclic, isomorphic to
$C_k$ for some $k$.  Thus
$G$ can be expressed as an extension
$$1\to C_k\to G\to C_q\to 1$$  Moreover a generator of $C_q$ acts on $C_k$ by
conjugation in $G$ so that $x\to x^r$ for some $r$ prime to $k$ such that $r^q\equiv
1\bmod k$.  In particular, $G$ has a presentation of the form
$$
\left<a,b:a^k, b^q=a^s, bab^{-1}=a^r\right>
$$ 
From this it follows that there are the following four possibilities:

\begin{enumerate}
\item $G$ is cyclic of order $kq$
\item $G$ contains a subgroup that is a nonabelian metacyclic group of order $pq$ with $p$ and $q$ distinct
primes
\item $G$ contains an elementary abelian subgroup of the form $C_p\times C_p$
\item $G$ contains a (generalized) quaternion $2$-group
\end{enumerate}

Possibilities 2 and 3 are ruled out by Propositions~\ref{metacyclic} and \ref{elemabel}, respectively.  Possibility 4 is ruled out
because a generalized quaternion group contains a nontrivial center, implying that all maximal subgroups intersect
nontrivially, a situation which was ruled out in Case 1.

We are left with the possibility that $G$ must be cyclic.  As before, the Lefschetz Fixed Point Formula implies that $G$ acts
semifreely, completing the proof.
\end{proof}

For completeness we explicitly bring the various pieces together to prove the Main Theorem.

\begin{proof}[Proof of the Main Theorem]
Suppose that the finite group $G$ acts pseudofreely, homologically trivially, locally linearly, preserving orientation, on a closed,
simply connected $4$-manifold $X$.  According to Theorem \ref{semifree} the action must be semifree.  But then, by Theorem
\ref{cyclic}, $G$ must be cyclic.  Finally, for any non-identity element $g\in G$ we have
$$
|X^G|=|X^g|=\chi(X^g)=\Lambda(g)=\chi(X)=b_2(X)+2
$$ 
\end{proof}

\section{\large Concluding Remarks}
One naturally wonders about removing the hypothesis of pseudofreeness in the Main Theorem. 
Since the $2$-torus $T^2$ acts on standard $4$-manifolds, such as connected sums of $\pm
\mathbf{C}P^2$ and $S^2\times S^2$ (see P. Orlik and F. Raymond \cite{Orlik70}), we see that
such manifolds do admit homologically trivial actions of $C_p\times C_p$, for example.  Of
course, the latter actions are definitely not semifree.  We conjecture that these are the only
sorts of groups that can act homologically trivially, when $b_2\ge 3$.  More precisely:

\begin{conj}
A finite group that acts locally linearly and
homologically trivially, on a closed, simply connected $4$-manifold $X$, with
$b_2(X)\ge 3$,  must be abelian of rank at most $2$.
\end{conj}

\bibliographystyle{amsplain}

\begin{thebibliography}{10}

\bibitem{Assadi90} Amir H. Assadi, Integral representations of finite
transformation groups III:  simply-connected four-manifolds, \emph{J. Pure Appl. Alg} {\bf
65} (1990), 209-233. 

\bibitem{Borel60}Armand Borel, \emph{Seminar on Transformation Groups,} Annals of
Math. Studies 46, Princeton University Press, 1960.
%

\bibitem{Bredon72}Glen E. Bredon, \emph{Introduction to Compact Transformation Groups,}
Academic Press, 1972.

\bibitem{Brown82}Kenneth S. Brown, \emph{Cohomology of Groups,} Springer-Verlag, New
York, 1982.

\bibitem{tomDieck87}Tammo tom Dieck, \emph{Transformation Groups,} Walter de Gruyter
\& Co., Berlin, 1987.

\bibitem{Edmonds87}Allan L. Edmonds, Construction of group actions on
four-manifolds, \emph{Trans. Amer. Math. Soc.} {\bf 299} (1987), 155-177.



\bibitem{Hambleton88} Ian Hambleton and Ronnie Lee, Finite group actions on
$P^2(\mathbf{C})$, \emph{ J. Algebra} {\bf 116} (1988), 227-242.



\bibitem{McCooey98}Michael McCooey, Groups that act pseudofreely on $S^2\times
S^2$, Indiana University preprint, 1998.

\bibitem{Milnor57} John Milnor, Groups which act on $S^n$ without fixed points, \emph{Amer. J. Math.} {\bf 79}
(1957), 623-630.

\bibitem{Orlik70} Peter Orlik and Frank Raymond, Actions of the torus on
$4$-manifolds I, \emph{Trans. Amer. Math. Soc.} {\bf 152} (1970), 531-559.

\bibitem{Darek87}Dariusz Wilczy\'{n}ski, Group actions on the complex
projective plane, \emph{Trans. Amer. Math. Soc. }{\bf 303}, (1987), 707-731.

\bibitem{Darek90}Dariusz Wilczy\'{n}ski, Symmetries of homology complex
projective planes, \emph{Math. Zeit. }{\bf 203}, (1990), 309-319.

\end{thebibliography}

\end{document}